\documentclass[leqno]{article}

\begin{document}

\title{Potpourri, 7}

\author{Stephen William Semmes	\\
	Rice University		\\
	Houston, Texas}

\date{}

\maketitle

	Let $V$ be a finite-dimensional real or complex vector space,
with positive dimension, and let $\|v\|_V$ be a norm on $V$.  Thus
$\|v\|_V$ is a nonnegative real number for all $v \in V$, $\|v\|_V =
0$ if and only if $v = 0$, $\|\alpha \, v\|_V = |\alpha| \, \|v\|_V$
for all real or complex numbers $\alpha$, as appropriate, and all $v
\in V$, and $\|v + w \|_V \le \|v\|_V + \|w\|_V$ for all $v, w \in V$.

	Associated to this norm we get a metric $\|v - w\|_V$ on $V$.
Since $V$ has finite dimension, it is linearly isomorphic to ${\bf
R}^n$ or ${\bf C}^n$, where $n$ is the dimension of $V$, according to
whether $V$ is a real or complex vector space.  The topology on $V$
determined by the metric $\|v - w\|_V$ corresponds exactly to the usual
topology on ${\bf R}^n$ or ${\bf C}^n$ with respect to this
isomorphism, by standard arguments.

	Let $\lambda$ be a linear functional on $V$, which is to say a
linear mapping from $V$ into the real or complex numbers, as
appropriate.  Using the identification with a Euclidean space, one can
see that there is a nonnegative real number $k$ such that
\begin{equation}
	|\lambda(v)| \le k \, \|v\|_V
\end{equation}
for all $v \in V$.

	Put
\begin{equation}
	\|\lambda\|_{V^*} = \sup \{|\lambda(v)| : v \in V, \  \|v\|_V \le 1\}.
\end{equation}
The supremum is finite by the remarks in the preceding paragraph, and
$\|\lambda\|_{V^*}$ can also be described as the smallest nonnegative
real number $k$ such that $|\lambda(v)|$ is less than or equal to $k$
times $\|v\|_V$ for all $v \in V$.  The space of all linear
functionals on $V$ is denoted $V^*$ and is a vector space with the
same dimension as $V$, and $\|\lambda\|_{V^*}$ defines a norm on $V^*$
called the dual norm associated to $\|v\|_V$.

	If $v$ is any vector in $V$, then a famous theorem implies
that there is a linear functional $\lambda$ on $V$ such that
$\|\lambda\|_{V^*} = 1$ and
\begin{equation}
	\lambda(v) = \|v\|_V.
\end{equation}
It follows that
\begin{equation}
	\|v\|_V = 
	   \max \{|\lambda(v)| : \lambda \in V^*, \  \|\lambda\|_{V^*} \le 1\}
\end{equation}
for all $v \in V$.

	Let $\Sigma$ denote the unit sphere in $V^*$ with respect to
the dual norm $\|\lambda\|_{V^*}$, i.e.,
\begin{equation}
	\Sigma = \{\lambda \in V^* : \|\lambda\|_{V^*} = 1\}.
\end{equation}
This is a compact subset of $V^*$ with respect to the topology induced
by the dual norm, which is the same as a topology induced from the
standard Euclidean topology using an isomorphism with Euclidean space.
Let us write $C(\Sigma)$ for the vector space of continuous real or
complex-valued functions on $\Sigma$, using the same scalars as for
$V$.

	For $\phi \in C(\Sigma)$ put
\begin{equation}
	\|\phi\| = \sup \{|\phi(\lambda)| : \lambda \in \Sigma\},
\end{equation}
where in fact the maximum is attained by standard results in analysis.
One can check that this defines a norm on $C(\Sigma)$, called the
supremum norm.

	For each $v \in V$ we get a function $\phi_v \in C(\Sigma)$
given by
\begin{equation}
	\phi_v(\lambda) = \lambda(v).
\end{equation}
The correspondence $v \mapsto \phi_v$ defines a linear mapping from $V$
into $C(V)$ such that
\begin{equation}
	\|\phi_v\| = \|v\|_V
\end{equation}
for all $v \in V$.

	Let $E$ be a finite nonempty set, and let $\mathcal{F}(E, V)$
denote the vector space of $V$-valued functions on $E$.  If $f$ is a
real or complex-valued function on $E$, put
\begin{equation}
	\|f\|_\infty = \max \{|f(x)| : x \in E \},
\end{equation}
which defines a norm on the vector space of real or complex-valued
functions on $E$.  Similarly, if $f$ is a $V$-valued function on $E$,
then put
\begin{equation}
	\|f\|_{\infty, V} = \max \{\|f(x)\|_V : x \in E\},
\end{equation}
which defines a norm on $\mathcal{F}(E, V)$.

	Suppose that $t(x, y)$ is a real or complex-valued function on
the Cartesian product $E \times E$.  Thus we get a linear mapping $T$
from the vector space of real or complex-valued functions on $E$ to
itself, as appropriate, given by
\begin{equation}
	T(f)(x) = \sum_{y \in E} t(x, y) \, f(y).
\end{equation}
Let $a$ be a nonnegative real number such that
\begin{equation}
	\|T(f)\|_\infty \le a \, \|f\|_\infty
\end{equation}
for all real or complex-valued functions $f$ on $E$.  This is
equivalent to the condition that
\begin{equation}
	\sum_{y \in E} |t(x, y)| \le a
\end{equation}
for all $x \in E_2$.

	Let us write $T_V$ for the version of $T$ acting on $V$-valued
functions on $E$, assuming that $V$ is real or complex according to
whether $t(x, y)$ is real or complex-valued.  Namely, for a $V$-valued
function $f$ on $E$, $T_V (f)(x)$ is defined by the same formula as in
the scalar case.  One can check that
\begin{equation}
	\|T_V (f)\|_{\infty, V} \le a \, \|f\|_{\infty, V}
\end{equation}
for all $V$-valued functions $f$ on $E$, using the condition on $|t(x,
y)|$ mentioned in the previous paragraph.  Alternatively, one can
consider the linear operator $T_{C(\Sigma)}$ which is the version of
$T$ on $C(\Sigma)$-valued functions on $E$, which basically acts
pointwise on $\Sigma$, and derive the analogous inequality using the
supremum norm on $C(\Sigma)$ from the inequality for scalar-valued
functions.  The inequality for $V$-valued functions then follows by
employing the isometric embedding of $V$ in $C(\Sigma)$.

	If $f$ is a real or complex-valued function on $E$, then put
\begin{equation}
	\|f\|_1 = \sum_{x \in E} |f(x)|,
\end{equation}
and more generally if $f$ is a $V$-valued function on $E$ put
\begin{equation}
	\|f\|_{1, V} = \sum_{x \in E} \|f(x)\|_V.
\end{equation}
These define norms on the vector spaces of scalar-valued and
$V$-valued functions on $E$.

	Suppose that $t(x, y)$, $T$ are as before, and that $b$ is a 
nonnegative real number such that
\begin{equation}
	\|T(f)\|_1 \le b \, \|f\|_1
\end{equation}
for all scalar-valued functions $f$ on $E$.  This is equivalent to the
condition that
\begin{equation}
	\sum_{x \in E} |t(x, y)| \le b
\end{equation}
for all $x \in E$.  In this event we have that
\begin{equation}
	\|T(f)\|_{1, v} \le b \, \|f\|_{1, V}
\end{equation}
for all $V$-valued functions $f$ on $E$ too.

	Fix a real number $p$, $1 < p < \infty$, and put
\begin{equation}
	\|f\|_p = \bigg(\sum_{x \in E} |f(x)|^p \bigg)^{1/p}
\end{equation}
for scalar-valued functions on $E$ and
\begin{equation}
	\|f\|_{p, V} = \bigg(\sum_{x \in E} \|f(x)\|_V^p \bigg)^{1/p}
\end{equation}
for $V$-valued functions on $E$.  Again these define norms on the
vector spaces of scalar and $V$-valued functions on $E$.

	Let us take for $V$ the vector space of real or complex-valued
functions on $E$, as appropriate.  A $V$-valued function on $E$ is
then basically the same as a real or complex-valued function $f(x, z)$
on $E \times E$.  The linear transformation $T_V$ on $V$-valued
functions on $E$ can be described explicitly by saying that if $f(w,
z)$ is a real or complex-valued function on $E \times E$, as
appropriate, then $T_V(f)$ is the function on $E \times E$ given by
\begin{equation}
	T_V(f)(x, z) = \sum_{y \in E} t(x, y) \, f(y, z).
\end{equation}

	Let us use the norm $\|\cdot \|_1$ on functions on $E$ as our
norm on $V$.  If $1 \le p < \infty$ and $f(x, z)$ is a function on $E
\times E$, which is equivalently a $V$-valued function on $E$, then
\begin{equation}
	\|f\|_{p, V} 
  = \bigg(\sum_{x \in E} \bigg(\sum_{z \in E} |f(x, z)| \bigg)^p \bigg)^{1/p}.
\end{equation}

	Suppose that $h(x)$ is a real or complex-valued function on
$E$, as appropriate.  Consider the function $H(x, z)$ on $E \times E$
given by $H(x, z) = h(x)$ when $x = z$, $H(x, z) = 0$ when $x \ne z$.
Thus we have that $\|H\|_{p, V} = \|h\|_p$ for all $p$, $1 \le p <
\infty$.

	Fix $p$, $1 < p < \infty$, and suppose that $B_p$ is a
nonnegative real number such that
\begin{equation}
	\|T_V(f)\|_{p, V} \le B_p \, \|f\|_{p, V}
\end{equation}
for all $V$-valued functions on $E$.  If $h$ is a scalar-valued
function on $E$ and $H$ is the associated function on $E \times E$
described in the previous paragraph, then $T_V(H)$ is the function on
$E \times E$ given by $t(x, z) \, h(z)$.  Hence
\begin{equation}
	\bigg(\sum_{x \in E} \bigg(\sum_{z \in E} |t(x, z)| \, |h(z)|
						\bigg)^p \bigg)^{1/p}
		\le B_p \, \bigg( \sum_{w \in E} |h(w)|^p \bigg)^{1/p}
\end{equation}
for all real or complex-valued functions $h$ on $E$.

\end{document}